\documentclass[11pt]{article}
\usepackage[utf8]{inputenc}
\usepackage{amsmath}
\usepackage{amsfonts}
\usepackage{amssymb}
\usepackage{mathrsfs}
\usepackage{tikz-cd}
\usepackage{amsthm}
\usepackage{mathtools}
\usepackage{hyperref}
\usepackage{graphicx}
\usepackage[labelfont=bf]{caption}
\usepackage{verbatim}
\usepackage{authblk}
\usepackage{amsmath, mathtools, amsfonts, mathrsfs, hyperref}

\graphicspath{ {./images/} }

\newcommand{\qp}{\mathbb{Q}_p}

\newcommand{\ap}{a_p}

\newcommand{\sT}{\mathcal{T}}

\newcommand{\til}[1]{\widetilde{#1}}

\newcommand{\br}[1]{\overline{#1}}

\newcommand{\sL}{\mathfrak{L}}
\newcommand{\brqp}{\br{\mathbb{Q}}_p}

\newtheorem{theorem}{Theorem}[section]

\newtheorem{conjecture}[theorem]{Conjecture}

\title{Zig-zag for Galois Representations}

\author{Eknath Ghate}
\affil{School of Mathematics, Tata Institute of Fundamental Research,\\ Homi Bhabha Road, Mumbai 400005, India}

\begin{document}
\maketitle
\begin{abstract}
  The zig-zag conjecture says that the reductions of two-dimensional crystalline representations
  of the Galois group of $\qp$ of large exceptional weights and half-integral slopes up to $\frac{p-1}{2}$
  vary through an alternating sequence of irreducible and reducible mod $p$ representations.
  We prove this conjecture in smoothly varying families of such representations
  for $p \geq 5$. The proof
  uses a limiting argument due to Chitrao-Ghate-Yasuda to
  reduce to the case of semi-stable representations of weights at most $p+1$, 
  and then appeals to the work of Breuil-M\'ezard, Guerberoff-Park and Chitrao-Ghate. 
\end{abstract}
  
\section{Introduction}

Let $p$ be a prime. There has been much interest in recent years on various
aspects of the $p$-adic and mod $p$ local Langlands program in higher dimensions
and over general local fields.
It is, however, remarkably still an open problem to describe the shape of the
reductions of
two-dimensional
crystalline representations $V_{k,a_p}$ of the Galois group $\mathrm{Gal}(\br{{\mathbb Q}}_p/\qp)$.
Here $a_p \in \brqp$ and $k \geq 2$ is an integer.
The resolution of this problem is important given that such representations
arise as the local $p$-adic Galois representations attached to modular forms with trivial
character, at primes not dividing the level of the form.

The first breakthroughs in describing the reduction
were made by Deligne, and Fontaine and Edixhoven. Their results can be stated in terms of the
{\it slope} $v_p(a_p)$ of $V_{k,a_p}$ where the $p$-adic valuation $v_p$ of $\brqp$ is
normalized so that $v_p(p) = 1$. The former showed that if $v_p(a_p) = 0$, then the reduction
of $V_{k,a_p}$ is reducible, whereas the latter pair showed that if $v_p(a_p) > 0$, then
the reduction is irreducible for the small weights $2 \leq k \leq p+1$. A well-known
conjecture due to Breuil, Buzzard, and Emerton says that if the weight $k$ is even and the
slope $v_p(a_p)$ is not an integer,
then the reduction of $V_{k,a_p}$ is always irreducible.

In this paper, we shall recall and prove (in families) our zig-zag conjecture which describes
the reduction of $V_{k,a_p}$ for large exceptional weights and half-integral slopes.
Let us assume from now on
that $p$ is an odd prime (eventually taken to be at least $5$).
Recall that a weight $k$ is said to be {\it exceptional} for the representation $V_{k,a_p}$ if
the slope $v_p(a_p)$ is a half-integer and $$k \equiv 2 v_p(a_p) + 2 \mod (p-1).$$
Exceptional weights fall squarely outside the scope of the conjecture in the
previous paragraph since such weights are even exactly when the slope $v_p(a_p)$ is
an integer.
In fact,
the zig-zag conjecture states that for large exceptional weights and half-integral
slopes in the range $$\frac{1}{2}  \leq v(a_p) \leq \frac{p-1}{2},$$ the
reduction of $V_{k,a_p}$ is given by an alternating sequence of irreducible and reducible
representations depending on the sizes of two auxiliary parameters $\tau$ and $t$.
In particular, it predicts the existence of several new cases where the
reduction is reducible.


A word on the provenance of the term `zig-zag'. Historically, when trying to compute
the reductions of $V_{k,a_p}$ using
the compatibility between the $p$-adic and mod $p$ local Langlands correspondences with
respect to the process of reduction, exceptional weights have turned out
to be the trickiest to treat.  In fact, the conjecture acquired its name because of
a zig-zag like pattern that emerges in the selection of Jordan-H\"older factors in
the associated graded of the theta filtration of the $(k-2)$-th symmetric power
representation of $\mathrm{GL}_2({\mathbb F}_p)$.

We recall the statement of the conjecture
\footnote{This version is mildly stronger than the original version 
   \cite[Conjecture 1.1]{Gha} in that there we require $k$ to be sufficiently far away
   from some weights which are strictly larger than $p+1$, whereas here
   $k$ is required to be sufficiently close to the weights $3 \leq k_0 \leq p+1$.}
 (see \cite{Gha}).

\begin{conjecture}[Zig-Zag Conjecture] 
  \label{zig-zag-conj}
  Say that 
  $k \equiv k_0 = 2v(a_p)+2 \mod (p-1)$ is a congruence class of exceptional weights
  for a particular half-integral slope $\frac{1}{2} \leq v_p(a_p) \in \frac{1}{2} {\mathbb Z} \leq \frac{p-1}{2}$. 
  Let $r = k-2$ and $r_0 = k_0-2$.
  Define two parameters:
  \begin{eqnarray*}
    \tau = v_p(c) & \text{and} &    t  = v_p(k - k_0),
  \end{eqnarray*}
  where
  \begin{eqnarray*}
    \label{c}
     c  & = &  \dfrac{a_p^2 - \binom{r-v_-}{v_+} \binom{r - v_+}{v_-} p^{r_0}}{p a_p}
 \end{eqnarray*}
 with $v_-$ and $v_+$ the largest and smallest integers 
 such that $v_p(a_p)$ lies in $(v_-, v_+)$.
  
Then, for all weights $k \geq k_0$  with $t$ sufficiently large,
the (semi-simplification of
 the) reduction $\bar{V}_{k,a_p}$ of the crystalline representation 
 $V_{k,a_p}$ on $\mathrm{Gal}(\brqp/\qp)$
 is given by:
 \begin{eqnarray*}
    \bar{V}_{k, a_p}  & \sim & 
  \begin{cases}
    \begin{array}{l}
      \mathrm{ind}(\omega_2^{r_0+1+i(p-1)}),    
    \end{array}                                                              & \text{ if } \tau \in (t + i-1, t + i)\\
    \begin{array}{l}
      \mu_{\lambda_i} \cdot \omega^{r_0-i}\,\oplus\, \mu_{\lambda_i^{-1}} \cdot \omega^{1+i},
    \end{array}                                                              & \text{ if } \tau = t  + i 
\end{cases}
\end{eqnarray*}
for $0 \leq i \leq \frac{r_0-1}{2}$ if $r_0$ is odd and $0 \leq i \leq \frac{r_0}{2}$ if $r_0$ is even and
for some constants $\lambda_i \in \br{{\mathbb F}}_p^*$.
In particular, as $\tau$ varies through the $t$-line,
the reduction $\bar{V}_{k,a_p}$ varies through an alternating sequence of irreducible and
reducible representations.
\end{conjecture}

We adopt the following conventions:
\begin{itemize}
  \item If $k > k_0$, so that $t$ is finite, then we interpret the
first interval  $(t-1,t)$ occurring above as $(-\infty,t)$.  Moreover, we interpret the last interval $[t+\frac{r_0 -1}{2}, t+\frac{r_0 -1}{2}]$ as $[t+\frac{r_0-1}{2},\infty]$
if $r_0$ is odd, and interpret the second last interval $(t + \frac{r_0-2}{2}, t+\frac{r_0}{2})$ as $(t + \frac{r_0-2}{2}, \infty]$ if 
$r_0$ is even and ignore the last interval $[t+\frac{r_0}{2}, t+\frac{r_0}{2}]$.
\item If $k = k_0$, so that $t = \infty$, then we interpret the first interval $(t-1,t)$ as $(-\infty,\infty]$ and ignore all other intervals.
\end{itemize}
Also, as usual, here $\omega$ is the mod $p$ cyclotomic character, `$\mathrm{ind}$' means
induction from $\mathrm{Gal}(\br{{\mathbb Q}}_p/{\mathbb Q}_{{p^2}})$
to $\mathrm{Gal}(\br{{\mathbb Q}}_p/\qp)$ and $\omega_2^c$ is the $c$-th power of the
fundamental character of level $2$ on the inertia subgroup $I_p$ extended to
$\mathrm{Gal}(\br{{\mathbb Q}}_p/{\mathbb Q}_{{p^2}})$  so that $\det(\mathrm{ind}(\omega_2^c)) = \omega^c$. Finally, $\mu_{\lambda_i}$ is the unramified
mod $p$ character taking geometric Frobenius to $\lambda_i$.

The zig-zag conjecture has so far only been proved\footnote{The first
  paper assumes $p \geq 3$ whereas the last two assume $p \geq 5$.
  The first work was done before the conjecture
was made. This and the second work provided the foundation on which to formulate the conjecture. The third paper
checked that the conjecture indeed holds in the next case. In fact, these papers establish a result stronger than Conjecture~\ref{zig-zag-conj}
in that there is no requirement 
that $t = v_p(k-k_0)$ be large (even $t = 0$ is allowed). However,
the assumption that $t$ is sufficiently large is necessary for larger slopes since there are counterexamples to this stronger version by Rozensztajn 
even for slope $2$ (see \cite{Gha}).} when the valuation $v_p(a_p)$ of 
$a_p$ is $\frac{1}{2}$, $1$ and $\frac{3}{2}$ \cite{BG13}, \cite{BGR18}, \cite{GR}.
The proof becomes more intricate as the slope increases, but in all cases the main
idea is to use the $p$-adic and mod $p$
local Langlands correspondences discovered by Breuil \cite{Br03}
along with a compatibility result of Berger \cite{B10} with respect to the process
of reduction, or more succinctly,
the functorial version of these correspondences due to Colmez \cite{Col10},
to reduce to computing the reduction
of the standard lattice in the completion of the locally algebraic representation
$\mathrm{ind}_{KZ}^G \mathrm{Sym}^{k-2}(E^2) / (T-a_p)$ of $G = \mathrm{GL}_2({\mathbb Q}_p)$.
Here $K$ is the compact open subgroup $\mathrm{GL}_2({\mathbb Z}_p)$ of $G$
and $Z ={\mathbb Q}_p^*$ is the center of $G$. This last computation
is rather tricky and involves using the Hecke operator $T$ to cut out the
above-mentioned Jordan-H\"older factors. 

We now turn to the main result of this paper.
Fix a `base' weight $3 \leq k_0 \leq p+1$. These weights vary through
a complete set of representatives modulo $(p-1)$.
We shall only consider weights $k$ which satisfy
$k \equiv k_0 \mod (p-1)$, that is, we work in
a fixed branch of weight space. We assume that the parameter
$$a_p = a_p(k)$$ appearing in $V_{k,a_p}$ is a {\it smooth} function
 of the weight $k$ (even $C^2$ in the sense of \cite{Sch84} suffices) in a $p$-adic
 neighborhood of $k_0$ and that $$a_p(k_0) = p^{\frac{k_0 -2}{2}}.$$
 In other words, we consider
a family $V_{k,a_p}$ of crystalline representations whose members are necessarily 
exceptional as $k \rightarrow k_0$ ($p$-adically) since the $p$-adic valuation of
$a_p(k)$ is eventually the same as that of $a_p(k_0)$.
The main result of this paper is the following theorem. It proves the
zig-zag conjecture for the members  of the family $V_{k,a_p}$
for $k$ sufficiently close to  $k_0$. 
 
 
\begin{theorem}
  \label{main-theorem}
  Let $p \geq 5$ and $3 \leq k_0 \leq p+1$. Suppose that $a_p = a_p(k)$ is a smooth function of
  $k \equiv k_0 \mod (p-1)$ in a $p$-adic neighborhood of $k_0$ and that 
  $a_p(k_0) = p^{\frac{k_0-2}{2}}$,
  so that $V_{k,a_p(k)} $ becomes exceptional as $k \rightarrow k_0$ of slope
  $\frac{1}{2} \leq v_p(a_p(k_0)) \leq \frac{p-1}{2}$.

  Then, the (semi-simplification of the) reduction $\bar{V}_{k,a_p(k)}$ of  $V_{k,a_p(k)}$
  on 
  the inertia subgroup $I_p$
  is exactly as described by 
  Conjecture~\ref{zig-zag-conj} for $k$ sufficiently close to $k_0$. Moreover,
  on $\mathrm{Gal}(\brqp/\qp)$, the constants $\lambda_i$ appearing in the unramified
  mod $p$ characters $\mu_{\lambda_i}$ can be determined
  completely (see \eqref{constants-main-theorem}, \eqref{star'}).
\end{theorem}

There are two main ingredients in the proof:
\begin{enumerate}
  \item [i)] Work of Chitrao-Ghate-Yasuda \cite{CGY} shows that the
 dual of $V_{k,a_p(k)}$ converges to the dual of the semi-stable
 representation $V_{k_0,\sL}$ for $\sL = 2 a_p'(k_0) / a_p(k_0)$ in a certain blow-up space of trianguline representations constructed by Colmez \cite{Col} and studied by Chenevier \cite{Che}, so $V_{k,a_p(k)}$ and $V_{k_0,\sL}$ have the same reduction as $k \rightarrow k_0$ with $k \neq k_0$.

\item [ii)] Recent work of Chitrao-Ghate \cite{CG23} computes
  the reduction of the semi-stable representations $V_{k_0,\sL}$
  for all weights
  $3 \leq k_0 \leq p+1$ and for all primes $p \geq 5$ revealing 
  a zig-zag like pattern in this setting. This 
  extends earlier important works of Breuil-M\'ezard \cite{BM}
  for even $k_0$ with $4 \leq k_0 < p+1$ on $\mathrm{Gal}(\br{{\mathbb Q}}_p/\qp)$ and  
  Guerberoff-Park \cite{GP} for odd $k_0$ with $3 \leq k_0 < p$ on $I_p$ (later
  extended partially to $\mathrm{Gal}(\br{{\mathbb Q}}_p/\qp)$ by Lee-Park \cite{LP}).\footnote{The works \cite{BM}, \cite{GP} were used in our first proof of zig-zag (https://arxiv.org/pdf/2211.12114.pdf) for slopes up to $\frac{p-3}{2}$.
    These works compute the reduction of the associated
    strongly divisible module which is not available for $k_0 > p$ (and
    partially so for weight $k_0 = p$).
  The paper \cite{CG23} computes instead the reduction
 of Breuil's Banach space $B(k_0,\sL)$ and then uses the Iwahori mod $p$ local Langlands
 correspondence \cite{Chi23}. In principle, this method works for all weights, and
 in particular, \cite{CG23} covers the cases $k_0 = p$ and $p+1$.}
\end{enumerate}
Putting i) and ii) together we immediately obtain a proof of zig-zag for the reduction of
$V_{k,a_p}$ for $a_p = a_p(k)$ for those weights $k$ that are  $p$-adically close to $k_0$. 

Thus, while \cite{CGY} used the zig-zag conjecture (and known cases of it for small slopes) to cross-check (and in some
cases of small weight to extend) the works \cite{BM}, \cite{GP}, the present paper is based on the observation that
one may use the limiting argument in \cite{CGY} to argue in the {\it reverse}  direction:
the works \cite{BM}, \cite{GP}, \cite{CG23} imply that zig-zag is eventually true for
the members of a family of crystalline representations because in the limit the
reductions behave like those of semi-stable representations. 

We now provide further details.

\section{Crystalline and semi-stable representations}

\subsection{Definitions}

Let $E$ be a finite extension of $\qp$. Let $D_{\mathrm{st}}$ be Fontaine's functor inducing an equivalence of categories
between semi-stable representations of the Galois group $\mathrm{Gal}(\br{{\mathbb Q}}_p/\qp)$ over $E$ and 
weakly admissible filtered $(\varphi, N)$-modules over $E$. The representations $V_{k,a_p}$ and $V_{k,\sL}$ are defined
using  this functor as follows.

For every integer $k \geq 2$ and $\ap \in E$ of positive valuation, $V_{k, \ap}$  is the irreducible two-dimensional {\it crystalline} representation over $E$ of the Galois group $\mathrm{Gal}(\br{{\mathbb Q}}_p/\qp)$ with Hodge-Tate weights $(0, k-1)$ such that $\mathrm{Fil}^i(D_{\mathrm{st}}(V^*_{k, a_p})) = E e_1$ for $1 \leq i \leq k - 1$, where $e_1$, $e_2$ are basis vectors over $E$ of the two-dimensional filtered $\varphi$-module $D_{\mathrm{st}}(V^*_{k, a_p})$, and such that $\varphi e_1   =p^{k-1}e_2$, $\varphi e_2 = -e_1 + a_p e_2$ (so that that the characteristic polynomial of $\varphi$ is $x^2 - a_px + p^{k - 1}$). 

Similarly, for every integer $k \geq 2$ and  $\sL \in E$ 
(called the $\sL$-invariant),   $V_{k, \sL}$ is the two-dimensional {\it semi-stable} representation over $E$ of the Galois group $\mathrm{Gal}(\br{{\mathbb Q}}_p/\qp)$  (which is known to be irreducible if $k \geq 3$) with Hodge-Tate weights $(0, k-1)$ and such that 
$\mathrm{Fil}^i(D_{\mathrm{st}}(V^*_{k, \sL})) = E(e_1 + \sL e_2)$ for $1 \leq i \leq k - 1$, 
where $e_1$, $e_2$ are basis vectors over $E$ of the two-dimensional filtered $(\varphi, N)$-module $D_{\mathrm{st}}(V^*_{k, \sL})$, and such that $\varphi e_1 = p^{k/2} e_1$, $\varphi e_2 = p^{(k-2)/2} e_2$ and $N e_1 = e_2$, $N e_2 = 0$. 


These representations  are important as they arise as 
the Galois representations $\rho_f |_{\mathrm{Gal}(\brqp/\qp)}$
attached to normalized cuspidal newforms $f$ of weight $k$ and trivial character:
\begin{itemize}
\item $V_{k,a_p}$ arises if $p$ is prime to the level of $f$ with $a_p$
      the $p$-th Fourier coefficient of the form
  \item $V_{k, \sL}$ arises if $p$ exactly divides the level of $f$ and the
        $p$-th Fourier coefficient of $f$ is $p^\frac{k-2}{2}$
        (as opposed to $- p^\frac{k-2}{2}$) with $\sL$ 
        the negative of the Fontaine-Mazur $\sL$-invariant of $f$.
\end{itemize}
\subsection{
  Limits of crystalline representations}

Let $E$ be a finite extension of $\qp$. The Robba ring $R$ over $E$ is the ring of bidirectional power series in $T$ over $E$ which converge on 
some clopen annulus $p^{-M} \leq |T| < 1$ ($M>0$ may vary with the power series).

Colmez \cite{Col} and Chenevier \cite{Che} have constructed a rigid analytic space whose $E$-valued points essentially parameterize
non-split trianguline $(\varphi, \Gamma)$-modules of rank $2$ over $R$. 
Let us recall the definition. Let $\sT$ be the rigid analytic variety of characters of $\qp^*$. Let $\chi$ be the $p$-adic cyclotomic
character and let $x : \qp^* \rightarrow \qp^*$ be the identity character. We call the
characters $x^i\chi$ for $i \geq 0$ {\it exceptional}.
For each $i \geq 0$, let $F_i$ and $F_i'$ be the closed analytic sub-varieties of $\sT \times \sT$
defined by 
\[\arraycolsep=2pt
    \begin{array}{rl}
        F_i(E) \! & = \{(\delta_1, \delta_2) \in \sT(E) \times \sT(E) \text{ } \vert \text{ } \delta_1\delta_2^{-1} = x^i\chi\}, \\
        F_i'(E) \! & = \{(\delta_1, \delta_2) \in \sT(E)\times\sT(E) \text{ } \vert \text{ } \delta_1\delta_2^{-1} = x^{-i}\}.
    \end{array}
  \]
Let $F = \cup_{i \geq 0} \> F_i$ and $F' = \cup_{i \geq 0} \> F_i'$. The Colmez-Chenevier space $\til{\sT}_2$ is the blow-up of $(\sT\times\sT) \setminus F'$ along $F$ in the category of rigid-analytic spaces. For background on blow-ups
 in the rigid setting
see \cite[\S 2]{CGY}.

The functor $D_{\mathrm{rig}}$ induces an equivalence of categories between representations of
$\mathrm{Gal}(\brqp/\qp)$ 
over $E$ and \'etale $(\varphi, \Gamma)$-modules over $R$ (see, e.g., \cite[\S 1]{Col}).
Moreover, $D_{\mathrm{rig}}(V_{k,a_p})$ and $D_{\mathrm{rig}}(V_{k,\sL})$ are  trianguline, so define points in $\til{\sT}_2$. The latter lie in the exceptional divisor of $\til{\sT}_2$.

The following theorem is proved in \cite[Remark 2 after Theorem 1.1]{CGY}.

\begin{theorem}
  \label{limit-theorem}
  Let $3 \leq k_0 \leq p+1$ and $k \equiv k_0 \mod (p-1)$. Assume that $a_p = a_p(k)$ is a smooth function of $k$ in a $p$-adic neighborhood of $k_0$
  and $a_p(k_0) = p^{\frac{k_0-2}{2}}$.

  Then, the sequence of crystalline representations $V^*_{k, a_p(k)}$ converges to the semi-stable 
  representation    $V^*_{k_0, \sL}$ in $\til{\sT}_2$ with $\sL = 2 a'_p(k_0)/a_p(k_0)$. In
  particular, for $k$ sufficiently close (but not equal) to $k_0$,
  the (semi-simplifications of the) reductions
  of $V_{k,a_p(k)}$ and $V_{k_0,\sL}$ are isomorphic. 
\end{theorem}

In the interest of keeping this paper self-contained, we sketch a proof of the theorem
in the rest of this section. For details see \cite{CGY}. The explicit triangulation of 
$D_{\mathrm{rig}}(V^*_{k,a_p})$ is 
\begin{eqnarray*}
  0 \rightarrow R(\delta_1(k)) \rightarrow D_{\mathrm{rig}}(V_{k,a_p(k)}^*) \rightarrow R(\delta_2(k)) \rightarrow 0
\end{eqnarray*}
where $\delta_1(k) = \mu_{y(k)}$ for $\mu_x$ the unramified character of $\qp^*$ taking $p$ to $x \in \brqp^*$,  
\begin{eqnarray}
  \label{formula-for-y(k)}
  y(k) = \frac{a_p(k) + \sqrt{a_p(k)^2  - 4p^{k-1}}}{2},
\end{eqnarray}
and $\delta_2(k) = \mu_{y(k)^{-1}}  \cdot \chi^{1-k}$, and where for a character $\delta$ of $\qp^*$, $R(\delta)$ is
the $(\varphi, \Gamma)$-module of rank 1 over $R$ with the twisted actions of $\varphi$ and $\Gamma$ given
by $\varphi f(T) = \delta(p) \cdot f((1+T)^p-1))$ and $\gamma f(T) = \delta(\chi(\gamma)) \cdot f((1+T)^{\chi(\gamma)}-1)$
for $f(T) \in R$ and $\gamma \in \Gamma$.

Since $a_p(k) \rightarrow a_p(k_0) = p^{r_0/2}$, formula \eqref{formula-for-y(k)} shows that $y(k)$ converges to $p^{r_0/2}$ as $k \rightarrow k_0$ with $k \neq k_0$. So 
$\delta_1(k) \rightarrow \delta_1$ and  $\delta_2(k) \rightarrow \delta_2$ for $\delta_1 :=  \mu_{p^{r_0/2}}$ and 
$\delta_2 :=  \mu_{p^{-r_0/2}} \chi^{1-k_0}$. Note that $\delta_1/\delta_2 = \mu_{p^{r_0}} \chi^{k_0-1} = x^{r_0} \chi$ is an exceptional character.
In particular, the sequence of crystalline representations $V^*_{k,a_p(k)}$ converges to a semi-stable representation
$V^*_{k_0,\sL}$ for some $\sL$ as $k \rightarrow k_0$. The computation of $\sL$ is
delicate and proceeds as follows.  

First it is not 
important to work with pairs of characters but only their quotients since the $\sL$-invariant only depends (by definition) on the quotient character.
Thus the first reduction is that we may work in $\tilde{\sT}$ which is the blow up of $\sT \setminus \{x^{-i} | i \geq 0\}$ along
the set $\{x^i \chi | i \geq 0\}$.

We pick a local affinoid chart $U_{r_0}$ in $\sT \setminus \{ x^{-i} | i \geq 0\}$ around $x^{r_0} \chi$. 
An explicit formula for the underlying affinoid algebra may be found in \cite[\S 2]{CGY}, but ignoring the tame part of our characters, the $E$-valued
points of $U_{r_0}$ are
\begin{eqnarray*}
  U_{r_0}(E) & = & \{ (x,y) \in p^{r_0} {\mathcal O}_E^* \times p {\mathcal O}_E \} \\
      \delta & \mapsto &  (\delta(p), \delta(1+p) - 1).
\end{eqnarray*}
Let $\tilde{U}_{r_0}$ be the preimage of $U_{r_0}$ in $\tilde{\sT}$. Again, explicit formulas for an admissible cover of 
$\tilde{U}_{r_0}$  may be found in {\it ibid.}, but the $E$-valued points of $\tilde{U}_{r_0}$ are
\begin{eqnarray*} 
  \tilde{U}_{r_0}(E) & = & \{ (x,y, a:b) \in U_{r_0}(E) \times {\mathbb P}^1(E) \> \vert \> (x-p^{r_0})b = (y - ((1+p)^{k_0-1}-1))a \}.
\end{eqnarray*}

Now the triangulation for the crystalline representations $V^*_{k,a_p(k)}$ above is unique (there is a unique non-split trianguline
extension of  $R(\delta_2(k))$ by $R(\delta_1(k))$ up to isomorphism) and in the coordinate patch for $\tilde{U}_{r_0}(E)$ 
above
it corresponds to 
\begin{eqnarray*}
  (y(k)^2, (1+p)^{k-1}-1, \> y(k)^2 - p^{r_0} : (1+p)^{k-1} - (1+p)^{k_0-1}) \in \tilde{U}_{r_0}(E).
\end{eqnarray*}
Dividing the third and fourth terms above by $k - k_0$ (for $k \neq k_0$), 
a short computation
shows that as $k \rightarrow k_0$ these points converge to a point in the exceptional divisor above $x^{r_0} \chi$,
namely:
\begin{eqnarray*}
  (p^{r_0}, (1+p)^{k_0 - 1} - 1, \> 2 a_p(k_0) a'_p(k_0): (1+p)^{k_0 - 1} \log(1+p))  \in \tilde{U}_{r_0}(E).
\end{eqnarray*}
Now a delicate computation shows that the $\sL$-invariant of the corresponding
$(\varphi, \Gamma)$-module in the sense\footnote{After modifying the sign.}  of Colmez \cite{Col} is $-\sL$ with
\begin{eqnarray*}
   \sL =  2\dfrac{a'_p(k_0)}{a_p(k_0)}.
\end{eqnarray*}
This computation involves working with two explicit bases of rank 2 trianguline $(\varphi, \Gamma)$-modules 
of exceptional weight, the first in terms of which the $\sL$-invariant is defined  \cite{Col},
and the second, considered by Benois \cite{Ben}, in terms of which the $\sL$-invariant can be computed.
We refer the reader to \cite[Theorem 5.2]{CGY} for the details.
In any case, the limiting Galois representation is $V(\delta_1, \delta_2, -\sL)$  in the notation of \cite{Col}  which
by comparing filtered $(\varphi, N)$-modules is isomorphic to $V^*_{k_0, \sL}$. This
proves the limiting statement in Theorem~\ref{limit-theorem}.

The final statement regarding the reductions follows from the following
three facts: \cite[Proposition 3.9]{Che} which shows that points
 in the rigid moduli
 space of trianguline representations $\til{\sT}_2$ live in families, a result of
 Kedlaya-Liu \cite{KL10} which
 shows that affinoid locally around an \'etale point such families come from families of Galois
 representations, and the discussion on \cite[p.~1513]{Che} which shows that 
 the semi-simplifications of the reductions of Galois representations living in connected families
 are isomorphic.

 This completes the overview of 
 the  proof of Theorem~\ref{limit-theorem}.
\qed

\section{Proof of zig-zag in the large}

We now prove the zig-zag conjecture for the representations $V_{k, a_p(k)}$ 
for $k$ sufficiently close to $k_0$ (Theorem~\ref{main-theorem}). For $k = k_0$, Conjecture~\ref{zig-zag-conj}
is anyway true by \cite{Edi} for all $3 \leq k_0 \leq p+1$. So we prove the theorem 
assuming $k \rightarrow k_0$ with $k \neq k_0$.

To this end, let us compute the parameter $\tau$ in the statement of Conjecture~\ref{zig-zag-conj}.
Recall that 
$v_{-}$ and $v_{+}$ are the largest and smallest 
integers such that $v_{-} < r_0/2 < v_{+}$. For $l \geq 1$, let $H_{l} = \sum_{i = 1}^{l}\frac{1}{i}$ be the
$l$-th partial harmonic sum, and set $H_0 = 0$. Set
$H_{-} = H_{v_-}$ 
and $H_{+} = H_{v_+}$. 
 For any $\sL$ in a finite extension of $\qp$, let
\[
    \nu = v_p(\sL - H_- - H_+)
\]
be the $p$-adic valuation of the $\sL$-invariant shifted by the partial harmonic sums $H_-$ and $H_+$.

Now since $a_p = a_p(k)$ is smooth about $k_0$, it has a Taylor expansion at $k_0$ (see, e.g., \cite[Proposition 28.3]{Sch84}),
which by Theorem~\ref{limit-theorem} can be written as
\begin{eqnarray*}
  a_p(k) & = & a_p(k_0) + a'_p(k_0)(k-k_0) + \cdots   
             \>\> = \> \>  p^{r_0/2} \cdot \left(1+ \frac{\sL}{2}(k-k_0) \right) + \cdots 
\end{eqnarray*}
where the $+ \cdots  $ means higher order terms involving at least $(k-k_0)^2$.
As $k \rightarrow k_0$, we have $v_p(a_p(k)) = r_0/2$. Thus for $k$ sufficiently close to $k_0$, we
have $v_{-} + v_{+} = k_0 - 2$. For $r = k - 2$, the first binomial coefficient in $\tau$ is given by 
    \begin{eqnarray*}
    {r - v_{-} \choose v_{+}} & = &   {k - k_0 + v_{+} \choose v_{+}} 
                                                 \> = \>  \frac{(k-k_0+v_+)(k-k_0+v_+-1))\cdots(k-k_0+1)}{v_{+}!} \nonumber\\
    & = & 1 +(k-k_0)\sum_{i = 1}^{v_{+}}\frac{1}{i} + \cdots  \>\> = \>\> 1 +(k-k_0) H_+ +  \cdots 
    \end{eqnarray*}
Similarly, the second binomial coefficient is
    ${r - v_{+} \choose v_{-}}  =  1 + (k-k_0) H_-  + \cdots$ 
Using the above formulas, we can compute  $\tau$ for $V_{k, a_p(k)}$ as $k \rightarrow k_0$ with $k \neq k_0$. We have
\begin{align}
  \label{tau-t-nu}
  \begin{split}  
    \tau & =  v_p\left(\frac{a_p(k)^2 - {r- v_{-}\choose v_{+}}{r - v_{+}\choose v_{-}}p^{r_0}}{pa_p(k)}\right) \\
           & =  v_p\left(\frac{p^{r_0}(1 + 2 \cdot \frac{\sL}{2}(k-k_0) + \cdots)  - p^{r_0}(1 +(k-k_0)(H_+ + H_- )+ \cdots )  }
                                      {p(p^{r_0/2}(1+\frac{\sL}{2}(k-k_0)) + \cdots)}\right) \\
           & =  r_0  + v_p((\sL - H_{-} - H_{+})(k-k_0) + \cdots )  - (1 + r_0/2) \\
           & =  r_0/2 - 1  + \nu  +  t,
 \end{split}
\end{align}
at least if $\nu$ is finite. When $\nu = \infty$, note that $\tau - t \rightarrow \infty$ as $k \rightarrow k_0$ with $k \neq k_0$.

But for $k$  sufficiently close to $k_0$ with $k \neq k_0$,  we have by Theorem~\ref{limit-theorem} that
$$\bar{V}_{k, a_p(k)}   \sim   \bar{V}_{k_0, \sL}$$
with $\sL = 2 a'_p(k_0) / a_p(k_0)$. 
Thus, for such $k$, we see that the reduction $\bar{V}_{k,a_p(k)} |_{I_p}$ is as described by Conjecture~\ref{zig-zag-conj} if and only if 
\begin{eqnarray*}
    \bar{V}_{k_0, \sL} |_{I_p} & \sim & 
  \begin{cases}
    \begin{array}{l}
      \mathrm{ind}(\omega_2^{r_0+1+i(p-1)}),    
    \end{array}                                                              & \text{ if } \nu \in (i-r_0/2, i-r_0/2+1)\\
    \begin{array}{l}
      \omega^{r_0-i}\,\oplus\,\omega^{1+i}, 
    \end{array}                                                              & \text{ if } \nu =  i -r_0/2+1,
\end{cases}
\end{eqnarray*}
for $0 \leq i \leq \frac{r_0-1}{2}$ if $r_0$ is odd and $0 \leq i \leq \frac{r_0}{2}$ if $r_0$ is even,
with conventions similar to those adopted after the statement of Conjecture~\ref{zig-zag-conj}. But this is
exactly what is proved by Breuil-M\'ezard \cite{BM} for even $r_0 \lneq p-1$, 
by Guerberoff-Park \cite{GP} for odd $r_0 \lneq p-2$  (cf. \cite[\S 1.4]{CGY}) and
by Chitrao-Ghate \cite[Theorem 1.1]{CG23} in the
respective missing cases $r_0 = p-1$ and $r_0 = p-2$. 
This proves zig-zag for the reductions of the crystalline representations
$V_{k,a_p(k)}$ 
for large weights $k$ sufficiently close to $k_0$ with $k \neq k_0$ on the inertia subgroup. The constants $\lambda_i$ appearing in the unramified characters in the reduction
on $\mathrm{Gal}(\brqp/\qp)$ are determined in the next section
  (see the formulas \eqref{constants-main-theorem}, \eqref{star'} below).
This proves Theorem~\ref{main-theorem}.
\qed

\section{Unramified constants}
  \label{section-constants}
  There remains the question of exact formulas for
  the constants $\lambda_i$ appearing in the unramified mod $p$ characters
  $\mu_{\lambda_i}$ in the reductions of the crystalline representations $V_{k,a_p}$
  in Conjecture~\ref{zig-zag-conj}. 

  In \cite{Gha}, we also conjectured that when $\tau = t+i$, 
  these constants have the shape:
\begin{align}
  \label{constants-crys}
  \begin{split}
  \lambda_i & =  \br{{*_i} \cdot \dfrac{c}{p^{i}}} 
                  \quad \text{ if } 0 \leq i < \dfrac{r_0 - 1}{2} \\
  \lambda_{i} + \lambda_i^{-1} & =  \br{{*_i}  \cdot \dfrac{c}{p^{i}}} 
  \quad \text{ if } i = \dfrac{r_0 - 1}{2} \text{ and } r_0 \text{ is odd},
  \end{split}
\end{align}
for some fudge factors $*_i$ which were left unspecified, but which were speculated  
to be simple rational expressions in $r$ with $r_0-r$ in the denominator. The 
following explicit formulas for the $*_i$ for the half-integral slopes
$v_p(a_p) = \frac{1}{2},1,\frac{3}{2}$ (so $r_0 = 1,2,3$, respectively)
were derived over the past decade (in \cite{BG13}, \cite{BGR18}, \cite{GR}, respectively)
by directly using the $p$-adic and mod $p$ local Langlands
correspondences in the crystalline setting:
\begin{eqnarray}
    \label{known-constants}
      *_0 = \frac{r_0}{r_0-r} \> \text{ for } r_0 = 1,2,3, \quad \text{ and }
    \quad  *_1 = \frac{r_0-1}{(r_0-r-1)(r_0-r)} \> \text{ for } r_0 = 3.
\end{eqnarray}

On the other hand, a large portion of \cite{CG23} is dedicated to computing
the corresponding constants $\lambda_i$ appearing in the reductions of the semi-stable
representations $V_{k_0, \sL}$ for all $1 \leq r_0 \leq p-1$ and $p \geq 5$,
including the tricky `self-dual' constant below.
When $\nu = i+1-r_0/2$, we have  by \cite[Theorem 1.1]{CG23},
\begin{align}
  \label{constants-st}
  \begin{split}
  \lambda_i & = \br{(-1)^{i+1} \> (i+1) {r_0-i \choose i+1}\dfrac{\sL - H_{-} - H_{+}}{p^{i+1-r_0/2}}},
               \quad       \text{ if } 0 \leq i < \dfrac{r_0 - 1}{2} \\
  \lambda_{i} + \lambda_i^{-1} & = 
       \br{(-1)^{i+1} \> (i+1){r_0 - i \choose i+1}\dfrac{\sL - H_{-} - H_{+}}{p^{i+1-r_0/2}}},
         \quad \text{ if } i = \dfrac{r_0 - 1}{2} \text{ and } r_0 \text{ is odd}.
       \end{split}
\end{align}
We can use \eqref{constants-st} to 
completely determine the constants $\lambda_i$ appearing in the reductions
of the family of crystalline representations $V_{k,a_k}$ with $a_p = a_p(k)$ in
Theorem~\ref{main-theorem}, for $k$ sufficiently close to $k_0$.
Indeed, the computation in \eqref{tau-t-nu} shows that for $\tau = t+i$
\begin{eqnarray*}
  \br{\dfrac{\sL - H_{-} - H_{+}}{p^{i+1-r_0/2}}}& = & \br{\dfrac{c}{p^i(r-r_0)}},
\end{eqnarray*}
as $k \rightarrow k_0$ with $k \neq k_0$. Thus, for all $1 \leq r_0 \leq p-1$,
\eqref{constants-st} yields the complete formulas:
\begin{align}
  \label{constants-main-theorem}
  \begin{split}
  \lambda_i & =  \br{ {*_i'}  \cdot \dfrac{c}{p^{i}}} 
                    \quad \text{ if } 0 \leq i < \dfrac{r_0 - 1}{2} \\
  \lambda_{i} + \lambda_i^{-1} & =  \br{{*_i'}  \cdot \dfrac{c}{p^{i}}} 
                                     \quad \text{ if } i = \dfrac{r_0 - 1}{2} \text{ and } r_0 \text{ is odd},
  \end{split}
\end{align}
with
\begin{eqnarray}
  \label{star'}
  *_i' & := &  \dfrac{ (-1)^{i} \> (i+1) {r_0-i \choose i+1}}{r_0-r}.
\end{eqnarray}
Comparing \eqref{star'} with the four known fudge factors $*_i$ in \eqref{known-constants}
one sees that
$$*_i' =  *_i,$$
in all but the last self-dual case, where the ratio of the two sides 
tends to $1$ as $k \rightarrow k_0$.

\vspace{.4cm}

{\bf \noindent Acknowledgements.} This work came about in view of a question asked by A. Pal after a talk
at a conference at ICTS in September 2022 organized by A. Burungale, H. Hida, S. Jha and Y. Tian.
A first proof for slopes up to $\frac{p-3}{2}$ was written up
in November 2022 during a visit to UNIST and KAIST with support from KAIX,
and I would like to thank C. Park and W. Kim for the invitations.
I also thank A. Chitrao for several
useful discussions over many years, and in particular, for recent joint
work which allowed us to extend the argument to the last two slopes 
$\frac{p-2}{2}$, $\frac{p-1}{2}$ in the conjecture.

\end{document}